\documentclass[12pt]{amsart}

\usepackage{txfonts}
\usepackage{amscd}
\usepackage{cite}

\usepackage{epsfig}
\usepackage{verbatim}
\usepackage{mathdots}
\usepackage{amssymb}
\usepackage{amsfonts}
\usepackage{amsbsy}
 \usepackage{graphicx}

\newtheorem{theo}{Theorem}[section]
\newtheorem{lem}[theo]{Lemma}
\newtheorem{prop}[theo]{Proposition}

\theoremstyle{definition}
\newtheorem{defi}[theo]{Definition}
\newtheorem{exam}[theo]{Example}
\newtheorem{rem}[theo]{Remark}

\numberwithin{equation}{section}

\newcommand{\R}{{\mathbb R}}
\newcommand{\Z}{{\mathbb Z}}
\newcommand{\CC}{{\mathbb C}}

\newcommand{\N}{{\mathbb N}}

\newcommand{\A}{{\mathcal A}}

\newcommand{\D}{{\mathcal D}}

\def\Re{\textrm{\textup{Re}}}
\def\Im{\textup{\textrm{Im}}}

\newcommand{\eproof}{\hfill$\square$}

\def\eqeproof{\eqno{\square}}

\begin{document}

\baselineskip 16pt

\title{ On the range of $\sum_{n=1}^\infty\pm c_n$}

\author[X.-G. He]{Xing-Gang He}
\address{School of Mathematics and Statistics   \\ Central China Normal University\\ Wuhan, 430079
\\P. R. China.}
\email{xingganghe@163.com; xingganghe@sina.com}
\author[C.-T. Liu]{Chun-Tai Liu}
\address{School of Mathematics and Statistics   \\ Central China Normal University\\ Wuhan, 430079
\\P. R. China.}
\email{lct984@163.com}
\thanks{
   \, \, This work was supported by the National Natural Science Foundation of China 10871180. 
}

\begin{abstract}
 Let $\{c_n\}_{n=1}^\infty$ be a sequence of complex numbers. In this paper we answer when the range of $\sum_{n=1}^\infty\pm
 c_n$ is dense or equal to the complex plane. Some examples are given to explain our results. As its
 application, we calculate the Hausdorff dimension of the level
sets of a Rademacher series with complex coefficients.
\end{abstract}

\maketitle

\maketitle \tableofcontents
\bigskip

\section{Introduction}
\setcounter{equation}{0}

Let $\{c_n\}_{n=1}^\infty$ be a sequence of complex numbers and let
$$Y_{\{c_n\}}=\sum_{n=1}^\infty \pm c_n,$$
where the `` + " and `` - " signs are chosen independently with
probability $1/2$. When all $c_n=a_n$ are real numbers, it is known
that $Y_{\{a_n\}}$ is a random variable if and only if  $\{a_n\}\in
\ell^2(\N)$, i.e., $\sum_{n=1}^\infty |a_n|^2<\infty$ \cite{jw}. In
this case, the distribution function of $Y_{\{a_n\}}$ is called the
{\it infinite Bernoulli convolution}, which has been studied
extensively from 1930's (see \cite{d,s} and the references given
there). It is clear that the support of the distribution function is
the whole real line if and only if $\{a_n\}\not\in \ell^1(\N)$. When
all $\{c_n\}\in \ell^2(\N)$ are complex numbers, $Y_{\{c_n\}}$ is
also a random variable. Clearly $\{c_n\}\not\in\ell^1(\N)$ does not
guarantee that
\begin{eqnarray}\label{main} R(\{c_n\}):=\left\{\sum_{n=1}^\infty \pm
c_n\right\}={\mathbb C}.
\end{eqnarray}
Motivating by this, in this paper we want to find rational
conditions such that \eqref{main} holds.

Another motivation for this issue is the Rademacher series (see
\cite{ast, dil, hu, mon, wx, xi}). A {\it complex Rademacher series}
associated to $\{c_n\}_{n=1}^\infty$ is defined by
$\sum_{n=1}^\infty c_n R(2^{n-1}x)$, where $R(x)$ is a periodic
function with period $1$ and $R(x)=\pm 1$ according to $x\in [0,
1/2)$ or $[1/2, 1)$, respectively. Clearly we have
\begin{eqnarray}\label{main1} R(\{c_n\})=\left\{\sum_{n=1}^\infty \pm
c_n\right\}=\left\{\sum_{n=1}^\infty c_n R(2^{n-1}x): x\in[0,
1)\right\}.
\end{eqnarray}

We cannot give a sufficient and necessary condition for the question
\eqref{main}.  In stead of it, we obtain a criterion for
$R(\{c_n\})$ being dense in the complex plane.
 \medskip

Let $c_n=a_n+ib_n\in\CC$ for $n\ge 1$ with $\{c_n\}=o(1)$, which
means that $\lim_{n\rightarrow \infty}c_n=0$. If $\{\alpha a_n+
\beta b_n
 \}_{n=1}^\infty\not\in
 {\ell}^1$
 for any  $\alpha, \beta\in \R$ with $\alpha+i\beta\ne 0$, we call the sequence $\{c_n\}_{n=1}^\infty$ a {\it linearly non-summable sequence}.

\begin{theo} \label{thm1}
    Let $\{c_n\}_{n=1}^\infty$ be a sequence in the complex plane. Then $R(\{c_n\})$ is dense in the complex plane $\CC$ if and only if
    $\{c_n\}_{n=1}^\infty$ is linearly non-summable.
\end{theo}

We are surprised that  there are some examples which satisfy
$\overline{R(\{c_n\})}=\CC$ but $R(\{c_n\})\ne \CC$ (Example
\ref{exam2}). At the same time, there are some examples with
$\overline{R(\{c_n\})}=\CC$ but we do not know whether they are
equal to $\CC$, an example with this property is $c_n=\frac
1{n\ln{n+1}}+\frac in$ for $n\ge 1$. The key step of the proof of
Theorem \ref{thm1} is the combination lemma (Lemma \ref{lem2.2}).

\medskip

 To give a sufficient condition for $R(\{c_n\})=\CC$ we
begin with a notation.
\begin{defi}
    Let $\{c_n=a_n+ib_n\}_{n=1}^\infty\not\in{\ell}^1$
    be a complex sequence with $\{c_n\}=o(1)$.  $t$ is called a {\it ratio} of $\{c_n\}_{n=1}^\infty$ if there exists a subsequence
    $\{c_{n_k}\}_{k=1}^\infty\not\in \ell^1$ such that $a_{n_k}/b_{n_k}\to t$ as $k\to \infty$, where $t$ may be infinity.
\end{defi}

\medskip

It is easy to check that a complex sequence $\{c_n\}_{n=1}^\infty$
is linearly non-summable if it has two distinct ratios.

\begin{theo}\label{thm2}
    Let $\{c_n\}_{n=1}^\infty$ be a sequence in the complex plane. Then $R(\{c_n\})$ is
    the complex space if $\{c_n\}_{n=1}^\infty$ has  two different ratios.
\end{theo}

 The difficult part of the proof of Theorem \ref{thm2} is how to show
 that $R(\{c_n\})$ contains an interior. We will use Moran function
 systems (Proposition \ref{prop3.1}) to overcome it.

\medskip

The other one interesting problem on this issue is to study the
level set of Rademacher series. As far back as 1930, Kaczmarz and
Steinhaus \cite{kac} showed that, for any $a\in\R$, the level set
    $$E_a:=\bigg\{x\in [0, 1): \sum_{n=1}^\infty a_n R(2^{n-1}x)=a\bigg\}$$
has continuous  cardinality if $\{a_n\}\not\in {\ell}^1$ and
$\{a_n\}=o(1)$. In 1962, Beyer \cite{bey} proved
    $$\dim_HE_a=1$$
under the assumption $\{a_n\}\in {\ell^2\setminus\ell^1}$. Wu
\cite{wu}  showed the same result under the conditions
$\{a_n\}\not\in {\ell}^1, \{a_n\}=o(1)$ and another man-made
condition $\sum_{n=1}^\infty|a_{n+1}-a_n|<\infty$. Finally, Xi
\cite{xi} obtained the result without the man-made condition.

In the complex case, we define the level set by $E_c=\{x\in [0, 1):
\sum_{n=1}^\infty c_n R(2^{n-1}x)=c\}$ for any $c\in\CC$.  we show
that
\begin{theo}\label{thm3}
    Let $\{c_n\}_{n=1}^\infty$ be a sequence in the complex plane.
    \begin{enumerate}
    \item If $\{c_n\}$ has two
        distinct
        ratios, then $\dim_HE_c=1$ for any $c\in \CC$.
    \item If $\{c_n\}$ is linearly non-summable with one ratio, then
        $$\dim_H\bigg\{x\in [0, 1): \sum_{n=1}^\infty c_n R(2^{n-1}x)\in B(c, \delta)\bigg\}=1$$
        for any $\delta>0$, where $B(c, \delta)$ is the ball with center at $c$ and radius $\delta$.
    \end{enumerate}
\end{theo}

\section{The combination lemma and proof of Theorem \ref{thm1}  }
\setcounter{equation}{0}

Let $\{c_n=a_n+ib_n\}_{n=1}^\infty\not\in {\ell}^1$ be a complex
sequence and let $\{-1, 1\}^\N$ be the set of  all sequences
$\{x_n\}_{n\in\N}$ satisfying $x_n\in\{-1, 1\}$ for $n\ge 1$. We
begin with the existence of ratios. Note that, in the definition of
a ratio, we demands that the subsequence is not in $\ell^1$.

\begin{prop}\label{lem2.1}
    Let $\{c_n=a_n+ib_n\}_{n=1}^\infty\not\in
   {\ell}^1$ be a  complex sequence with $\{c_n\}=o(1)$. Then there exists at least one ratio.
\end{prop}
\begin{proof}
Let $\Lambda_1=\{n: |\frac {a_n}{b_n}|\le 1\}$ and $\Lambda_2=\{n: |\frac {a_n}{b_n}|\ge 1\}$. Then, by the symmetry of $a_n$ and $b_n$, without
loss of generality we assume that $\sum_{n\in\Lambda_1}|c_n|=\infty$.

Let $\A_{j, \,k}=\{n: a_n/b_n\in [j/2^k, (j+1)/2^k]\}$ for $k\ge 0$
and $-2^k\le j<2^k$. Then, for each $k$, there exists $j_k$ such
that $\{[j_k/2^k, (j_k+1)/2^k]\}_{k=0}^\infty$ is a decreasing
sequence of sets and  $\sum_{n\in \A_{j_k, \, k}}|c_n|=\infty$ for
each $k\ge 0$.

We claim that $t_0$ is a ratio of $\{c_n\}_{n=1}^\infty$, where
$t_0=\lim_{k\to\infty}\frac{j_k}{2^k}$. We prove the claim as
follows: Choose a finite set ${\mathcal B}_1$ from $\A_{j_1, \, 1}$
so that $\sum_{n\in{\mathcal B}_1}|c_n|\ge1$. Then choose ${\mathcal
B}_2$ from $\A_{j_2, \, 2}$ so that $\min\{n: n\in {\mathcal B}_2\}>
\max\{n: n\in {\mathcal B}_1\}$ and $\sum_{n\in{\mathcal
B}_2}|c_n|\ge 1$. As so on we can choose $\{{\mathcal
B}_k\}_{k=1}^\infty$ satisfying $\min\{n: n\in {\mathcal B}_{k+1}\}>
\max\{n: n\in {\mathcal   B}_k\}$  and $\sum_{n\in{\mathcal
B}_k}|c_n|\ge 1$ for each $k\ge 1$. Then the sequence $\{c_n: n\in
\cup_{k=1}^\infty {\mathcal B}_k\}$ has the ratio $t_0$.
\end{proof}

The following result guarantees that one can replace   any two
different ratios with two preconcerted ratios by acting a linear
mapping.
\begin{prop}\label{prop2.1}
    Let $\{c_n=a_n+ib_n\}_{n=1}^\infty$ be a complex sequence and let
    $\begin{pmatrix}
        \alpha_1 & \beta_1\\
        \alpha_2 & \beta_2 \\
    \end{pmatrix}$
    be a non-singular matrix. Then $R(\{a_n+ib_n\})$ is dense in (equal to) $\CC$ if and only if $R(\{(\alpha_1a_n+\beta_1 b_n)+i(\alpha_2a_n+ \beta_2b_n)\})$ is dense in (
    equal to, resp.) $\CC$
\end{prop}
\proof
    The assertion follows from the identity:
    $$R(\{a_n+ib_n\})
    \begin{pmatrix}
        \alpha_1 & \beta_1\\
        \alpha_2 & \beta_2 \\
    \end{pmatrix}
    =R(\{(\alpha_1a_n+\beta_1 b_n)+i(\alpha_2a_n+\beta_2b_n)\}).
    \eqeproof$$

\bigskip

From now on, for any $c=a+ib\in \CC$, we use the norm
$\|c\|=\max\{|a|,|b|\}$ throughout  this paper. For any complex
sequence $\{c_n\}_{n\in I}$ where $I\subseteq \N$, we denote that
$$\|\{c_n\}_{n\in I}\|=\sup_{n\in I}\|c_n\|.$$

\medskip

The following fact will be used in the proof of Lemma \ref{lem2.2}:
Let $\|c_1=a_1+ib_1\|\le 1$ and $\|c_2=a_2+ib_2\|\le 1$. Then it is
easy to check that: {\it $\|c_1\pm c_2\|> 1$ is equivalent to that
$|a_1|+|a_2|> 1$, $|b_1|+|b_2|>1$ and $a_1a_2b_1b_2<0$.} The
following combination lemma plays  a key role in the proof of
Theorem \ref{thm1}.
\begin{lem} [combination lemma]\label{lem2.2}
    Let $\{c_n\}_{n=1}^5$ be complex numbers satisfying $\|c_n\|\le 1$ for $1\le n\le 5$ and $\|c_n\pm c_{n+1}\|> 1$ for $1\le n\le 4$.
    Then there exists $\{x_n\}_{n=1}^5\in \{-1, 1\}^5$ such that
        $$\bigg\|\sum_{n=1}^5x_nc_n\bigg\|\le 2.$$
\end{lem}
\def\sign{\textrm{\textup{sign}}}
\begin{proof}
    Let $\sign(x)$ be the sign function, that is, $\sign(x)=-1$, 0 and $1$ according to $x<0$, $x=0$ and $x>0$ respectively.
    Denote $u=c_1\sign(a_1)-c_2\sign(a_2)-c_3\sign(a_3)+c_4\sign(a_4)$. By the above fact and the hypotheses, we have the imaginary
    part of $u$ satisfying $|\Im(u)|=
  \big|
    |b_1|+|b_2|-|b_3|-|b_4|
   \big|\le 1$.
    Similarly, write $v=c_2\sign(a_2)-c_3\sign(a_3) -c_4\sign(a_4)+c_5\sign(a_5)$, we have $|\Im(v)|\le 1$.
\medskip

    We claim that, if the real part of $u$ satisfies $|\Re (u)|>1$, then $|\Re(v)|\leq1$. Since $\Re (u)=|a_1|-|a_2|-|a_3|+|a_4|$
    and $|a_2|+|a_3|>1$, the condition $|\Re (u)|>1$ implies that $|a_1|-|a_2|-|a_3|+|a_4|<-1$. Similarly, if $|\Re (v)|>1$, then
    $|a_2|-|a_3|-|a_4|+|a_5|<-1$. Consequently, $|a_1|-2|a_3|+|a_5|<-2$, which leads to $2<2|a_3|\le 2$ and it is impossible. Hence the claim follows.

    The result follows by choosing $\{x_n\}_{n=1}^5\in \{-1, 1\}^5$ such that $\sum_{n=1}^5x_nc_n=u+c_5$ when $|\Re (u)|\le 1$ or
    $\sum_{n=1}^5x_nc_n=c_1+v$ when $|\Re (u)|> 1$.
\end{proof}

\begin{lem}  \label{lem2.3}
    Let $\{c_n\}_{n=1}^N$ be complex numbers with all $\|c_n\|\le 1$. Then there
    exists
    $\{x_n\}_{n=1}^N\in \{-1, 1\}^N$ such that
        $$\bigg\|\sum_{n=1}^kx_nc_n\bigg\|\le 5, \qquad \mbox{for all}\, \, 1\le k\le N.$$
\end{lem}
\begin{proof}
    We prove it by induction. Assume that the result holds for $N$ and $N>5$. For $N+1$, we show it by  two cases:

    Case 1. There exists $j$, $1\le j\le4$, such that
        either $\|c_j+c_{j+1}\|\le 1$ or $\|c_j-c_{j+1}\|\le 1$. Without loss of generality we say that $\|c_1+c_{2}\|\le 1$. By induction there exists $\{x_i\}_{i=1}^N\in \{-1, 1\}^N$ such that
        $$\bigg\|x_1(c_1+c_2)+\sum_{i=2}^kx_ic_{i+1}\bigg\|\le 5$$
    for $2\le k\le N$, this implies the assertion for $N+1$;

 Case 2. The assumption in Case 1 fails. We replace
    $\{c_1, c_2, c_3, c_4\}$ by $u$ if $\|u\|\le 1$ or $\{c_2, c_3, c_4,c_5\}$ by $v$
   if $\|u\|> 1$, where $u$ and $v$ are given in the proof of Lemma \ref{lem2.2}.
    Then the assertion follows by Lemma \ref{lem2.2} and the same idea of Case 1.
\end{proof}

\bigskip

    Now we can give a result on the controlling problem.
\begin{prop}\label{coro2.1}
    Let $\{c_n\}_{n=1}^\infty$ be a sequence in the complex plane with $\{c_n\}=o(1)$. Then there exists a sequence $\{x_n\}_{n=1}^\infty\in \{-1, 1\}^\N$ such that
        $$\bigg\|\sum_{n=1}^\infty x_nc_n\bigg\|\le 5\|\{c_n\}_{n=1}^\infty\|.$$
\end{prop}
\proof
    Since $\{c_n\}=o(1)$ as $n\to \infty$, there exists an increasing natural number sequence $\{N_k\}_{k=0}^\infty$ such that $N_0=1$ and
        $$\bigg\|\{c_n\}_{n=N_{k-1}}^{N_{k}-1}\bigg\|\le 2^{-k}\|\{c_n\}_{n=1}^\infty\|$$
    for $k\ge 1$. Clearly the result follows by using Lemma \ref{lem2.3} for each subsequence $\{c_n\}_{n=N_{k-1}}^{N_{k}-1}$.
\eproof

\begin{lem}\label{lem2.6}
    Let $\{c_n=a_n+ib_n\}_{n=1}^\infty\subset \CC$ be linearly
    non-summable. If there exists a non-summable subsequence of
    $\{c_n\}_{n=1}^\infty$ such that its real or imagine part is
    summable,
    then $R(\{c_n\})$ is dense in $\CC$.
    .
\end{lem}
\begin{proof}
    Let $\{c_{n_k}\}_{k=1}^\infty$ be a  non-summable  subsequence
    satisfying that its real or imagine part is summable. Without loss of generality we assume that its real part $\sum_{k=1}^\infty |a_{n_k}|$
    converges and $n_1>1$.

    Note that $\{a_n\}_{n=1}^\infty$ is not in $\ell^1$ with $a_n\to 0$. For any $a+bi\in \CC$, there exists
    a sequence $\{x_n\}\in \{-1,1\}^\N$ such that $a=\sum_{n=1}^\infty a_nx_n$. We denote $B_k=\sum_{n=1}^{n_k-1} b_nx_n$ for $k\ge 1$,
    where $\{n_k\}$ is given in the above subsequence. Let  $\Lambda_k=\{n_k, n_{k+1}, n_{k+2}, \ldots\}$ and $\Lambda_k^c=\N\setminus
    (\Lambda_k\cup \{1, 2, \ldots, n_k-1\})$.

    Since $\{b_{n\in \Lambda_k}\}$ is not summable, there exist $y_n\in\{-1, 1\}$ for $n\ge n_k$ dependent on $k$  such that
     $\sum_{n\in\Lambda_k}b_ny_n=b-B_k$ and $\|\sum_{n\in\Lambda_k^c}c_ny_n\|\le 5\|\{c_n\}_{n\in\Lambda_k^c}\|$. Hence we have
    \begin{eqnarray*}
        \sum_{n=1}^{n_k-1}c_nx_n+\sum_{n=n_k}^\infty c_ny_n&=&a-\sum_{n=n_k}^\infty
        a_nx_n+iB_k+\sum_{n\in\Lambda_k}a_ny_n +i(b-B_k)+\sum_{n\in\Lambda^c_k}c_ny_n\\
        &=&a+bi-\sum_{n=n_k}^\infty
        a_nx_n+\sum_{n\in\Lambda_k}a_ny_n +\sum_{n\in\Lambda^c_k}c_ny_n.
    \end{eqnarray*}
    Note that
        $$\bigg|-\sum_{n=n_k}^\infty a_nx_n+\sum_{n\in\Lambda_k}a_ny_n+ \sum_{n\in\Lambda^c_k}c_ny_n\bigg|\le
        \bigg|\sum_{n=n_k}^\infty a_nx_n\bigg|+\sum_{n\in\Lambda_k}
        |a_n|+5\|\{c_n\}_{n\in\Lambda_k^c}\|, $$
    which tends to zero when $k$ tends to infinity. Then the proof  is complete.
\end{proof}

\bigskip

{\noindent \bf Proof of Theorem \ref{thm1}.}
    We first prove the sufficiency. Suppose that $\{c_n\}_{n=1}^\infty$ is linearly non-summable. By Lemma \ref{lem2.1} and Proposition \ref{prop2.1}
    there exists a subsequence $\{c_{n_k}\}_{k=1}^\infty\not\in\ell^1$ such that $a_{n_k}/b_{n_k}$ tends to $0$ when $k$ tends to infinity.
    In this case we still have $\sum_{n=1}^\infty |a_n|=\infty$ by the linear non-summation. When $\sum_{k=1}^\infty |a_{n_k}|$ converges,
    the sufficient condition follows by Lemma \ref{lem2.6}; When $\sum_{k=1}^\infty |a_{n_k}|$ diverges, we construct a subsequence
    $\{l_k\}$ of $\{n_k\}$ such that $\sum_{k=1}^\infty|a_{l_k}|<\infty$ and $\sum_{k=1}^\infty |b_{l_k}|=\infty$. This implies the sufficiency
    according to Lemma \ref{lem2.6} again.

   Now we construct a desired subsequence of $\{c_{n_k}\}_{k=1}^\infty$ if $\sum_{k=1}^\infty |a_{n_k}|$
   diverges.  Note that in this case $\sum_{k=1}^\infty |b_{n_k}|$ diverges also.
    Denote $\Lambda_m=\{k: |a_{n_k}/b_{n_k}|<2^{-m}\}$ for $m\ge 1$. Then $k$ belongs to $\Lambda_m$ for sufficiently large $k$ and thus
   $\sum_{k\in \Lambda_m} |b_k|=\infty$ for each $m$. We can choose
     finite sets $\Gamma_k\subset \Lambda_k$ such that $\sum_{n\in\Gamma_k}|b_n|\in(1, 2)$ and $\max\Gamma_k<\min\Gamma_{k+1}$ for $k\ge 1$.
      We claim that $\Gamma:=\cup_{k=1}^\infty \Gamma_k$ is the index of a desired subsequence. In fact,
        $$\sum_{n\in\Gamma}|a_n|=\sum_{k=1}^\infty \sum_{n\in\Gamma_k}|a_n|\le \sum_{k=1}^\infty 2^{-k}\sum_{n\in\Gamma_k}|b_n|\le 2$$
    and
        $$\sum_{n\in\Gamma}|b_n|=\sum_{k=1}^\infty \sum_{n\in\Gamma_k}|b_n|\ge \sum_{k=1}^\infty 1=\infty.$$

    Now we prove the necessity. Suppose that $R(\{c_n\})$ is dense in $\CC$, then  both the real and imaginary parts of $\{c_n\}$ are not in $
   {\ell}^1$.
    If there exist $\alpha$ and $\beta$ such that $\{\alpha a_n+\beta b_n\}_{n=1}^\infty \in
    {\ell}^1$,
    this implies a contradiction by Proposition \ref{prop2.1}.
\eproof

\section{Moran function systems and proof of Theorem \ref{thm2}}
\setcounter{equation}{0}

In the proof of Theorem \ref{thm2}, the difficult point is to show
that $R(\{c_n\})$ contains interiors. The following proposition
 will help us to show it \cite{feng}. We begin with
some notations.

  Given a natural number sequence
$\{n_k\}_{k=1}^\infty$ with all $n_k\ge 2$ and a sequence $\{f_{k,\,
i}(x):  k\ge 1, i=1, 2, \ldots, n_k\}$ of functions from $\R^n$ to
itself, which satisfy that
    $$\|f_{k,\, i}(x)-f_{k, \,i}(y)\|\le r\|x-y\|$$
for all $k\ge 1$ and $1\le i\le n_k$, where $0<r<1$. We say the
sequence a {\it Moran function system with contraction $r$}. Define
$\prod_{k=1}^m\{1,2,\ldots,n_k\}
=\{\sigma=\sigma_1\sigma_2\cdots\sigma_m: \sigma_k\in\{1, 2,\ldots,
n_k\},\,1\le k\le m\}$ for $m\ge 1$ and
$\prod_{k=1}^\infty\{1,2,\ldots,n_k\}
=\{\sigma=\sigma_1\sigma_2\cdots\cdots: \mbox{each}\, \,
\sigma_k\in\{1,2,\ldots,n_k\}\}$. For any
$\sigma=\sigma_1\cdots\sigma_m$, we define
   $$f_\sigma(x)=f_{1,\,\sigma_1}\circ f_{2,\,
   \sigma_2}\circ\cdots {\circ}f_{m, \, \sigma_m}(x), $$
which is the composing function of $f_{i,\,\sigma_i},i=1,2,\ldots, m$.
\begin{prop} \label{prop3.1}
    Let $F=\{f_{k,i}(x): i=1,2,\ldots,n_k,k\ge1\}$ be a Moran function system with contraction $r$. Suppose that the set $\{f_{k,i}(0): i=1,2,\ldots,n_k,k\ge1\}$ is bounded with bound $M$, then
    \begin{enumerate}
    \item For any $\sigma=\sigma_1\sigma_2\cdots\in \prod_{k=1}^\infty\{1,2,\ldots,n_k\}$, the limit
            $$\lim_{k\to \infty}f_{\sigma_1\cdots\sigma_k}(0)$$
        exists, we denote the value as $f_\sigma(0)$;
    \item The set $K_F$, $K_F:=\{f_\sigma(0): \sigma\in \prod_{k=1}^\infty\{1,2,\ldots,n_k\}\}$, is a nonempty compact set;
   \item Let $Q$ be a compact set so that $Q\subseteq \bigcup_{i=1}^{n_k}f_{k,i}(Q)$ for all $k\ge1$, then $Q\subseteq K_F$.
    \end{enumerate}
 \end{prop}
 \begin{proof}
    Let $B(0, R)$ be the closure ball with center at $0$ and radio $R$. Then, for each $k$ and $i$ with $1\le i\le n_k$, we have
        $$f_{k,\,i}(B(0,R))\subseteq B(f_{k,\,i}(0),rR) \subseteq B(0,M+rR)\subseteq B(0,R)$$
    if $R>M/(1-r)$. This implies that
    $\{\bigcup_{\sigma\in \prod_{i=1}^m\{1,2,\ldots,n_i\}} f_\sigma(B(0,R))\}_{m=1}^\infty$ is a decreasing sequence of compact sets. Hence,
        $$\bigcap_{m=1}^\infty\bigcup_{\sigma\in \prod_{i=1}^m\{1,2,\ldots,n_i\}}f_\sigma(B(0,R))$$
    is a nonempty compact set, which is independent of large $R$. For any $\sigma=\sigma_1\sigma_2\cdots\in \prod_{k=1}^\infty\{1,2,\ldots,n_k\}$, we have
        $$\bigcap_{k=1}^\infty f_{\sigma_1\cdots\sigma_k}(B(0,R))=
       \Big\{ \lim_{k\to\infty}f_{\sigma_1\cdots\sigma_k}(0)=f_\sigma(0) \Big\}.$$
    Hence,
    \begin{equation}\label{ine2}
        K_F=\bigcap_{m=1}^\infty \bigcup_{\sigma\in \prod_{i=1}^m\{1,2,\ldots,n_i\}}f_\sigma(B(0,R)).
    \end{equation}
    This deduces (1) and (2).

    By the hypothesis in (3) we have
        $$Q\subseteq \bigcup_{\sigma\in\prod_{k=1}^m\{1,2, \ldots,n_k\}}f_\sigma(Q)$$
    for $m\ge 1$. Choosing $R$ so that $Q\subseteq B(0, R)$, we obtain $Q\subseteq K_F$ by the above and \eqref{ine2}.
\end{proof}
\bigskip

Next we construct a Moran function system $F$ and a cube
$Q=[-5,5]\times [-5, 5]$ such that $Q\subseteq K_F\subseteq
R(\{c_n\})$ if $\{c_n\}_{n=1}^\infty$ has  two distinct ratios.

\medskip

The next lemma says that we can construct a subsequence from
$\{c_n\}$ such that its real part and imaginary part are comparable
to the sequence $\{\delta^n\}_{n=1}^\infty$, where $0<\delta<1$.

\begin{lem}\label{lem3.1}
    Let $\{c_n=a_n+ib_n\}_{n=1}^\infty\not\in
   {\ell}^1$  be a complex sequence with $\{c_n\}=o(1)$ and $\lim_{n\to\infty}a_n/b_n=t$, where $0<t<\infty$. Then, for any $0<\delta<1$ and a
   positive number
    sequence $\{\eta_k\}_{k=1}^\infty$ with $\sum_{k=1}^\infty\eta_k<\infty$, there exists an increasing sequence of sets $\{\Lambda_k\}_{k=1}^\infty$
     such that
            $$\bigg|\sum_{n\in\Lambda_k}|b_n|-\delta^k\bigg|\le \eta_k, \quad \bigg| \frac{\sum_{n\in\Lambda_k}a_n\sign(b_n)}
             {\sum_{n\in\Lambda_k}|b_n|}-t\bigg|<\eta_k$$
    for all $k\ge 1$ and
        $$\sum_{k=1}^\infty\sum_{n\in\Lambda_k}\|c_n\|<\infty.$$
\end{lem}
\proof
    We define $\Gamma_k=\{n: |a_n/b_n-t|+\|c_n\|<\eta_k\}$ for $k\ge 1$. Then by hypotheses we have
        $$\sum_{n\in\Gamma_k}|b_n|=\infty.$$
    Now we choose finite sets $\Lambda_k$ from $\Gamma_k$ by induction. Note that $b_n$ tends to zero when $n$ tends to infinity, we can choose a
    finite set $\Lambda_1$ from $\Gamma_1$ such that $\big|\sum_{n\in\Lambda_1}|b_n|-\delta\big|\le \eta_1, $ and then choose a finite set $\Lambda_2$
     from $\Gamma_2$ such that $\min\Lambda_2>\max\Lambda_1$ and $\big|\sum_{n\in\Lambda_2}|b_n|-\delta^2\big|\le \eta_2.$ In general we
     obtain an increasing sequence $\{\Lambda_k\}_{k=1}^\infty$ of finite
     sets and
        $$\big|\sum_{n\in\Lambda_k}|b_n|-\delta^k\big|\le \eta_k, \qquad \mbox{for}\, \, k\ge 1.$$

    Now we show the second inequality. Write $t_n=a_n/b_n$, then $|t_n-t|<\eta_k$ for $n\in\Gamma_k$. Hence, the second assertion follows from that
        $$\frac{\sum_{n\in\Lambda_k}a_n\sign(b_n)}{\sum_{n\in\Lambda_k}|b_n|}-t=\sum_{n\in\Lambda_k}(t_n-t)\frac{|b_n|}{\sum_{m\in\Lambda_k}|b_m|}.$$
    Note that
        $$\bigg|\sum_{n\in\Lambda_k}a_n\sign(b_n)\bigg|\le            (|t|+\eta_k)\sum_{n\in\Lambda_k}|b_n|\le (|t|+\eta_k)(\delta^k+\eta_k), $$
    then
        $$\sum_{k=1}^\infty\sum_{n\in\Lambda_k}\|c_n\|\le (1+|t|+\sum_{k=1}^\infty \eta_k)\sum_{k=1}^\infty (\delta^k+\eta_k)<\infty.
        \eqeproof$$

Applying Lemma \ref{lem3.1} to a special case, we have the following
lemma.

\begin{lem} \label{lem3.2}
    Let $\{c_n=a_n+ib_n\}_{n=1}^\infty\not\in
   {\ell}^1$ and $\{\gamma_n=\alpha_n+i\beta_n\}_{n=1}^\infty\not\in
   {\ell}^1$ with
    $\{c_n\},\{\gamma_n\}=o(1)$.
    Suppose $\lim_{n\to\infty}a_n/b_n=2$ and $\lim_{n\to \infty}\beta_n/\alpha_n=3$, then there exist two increase sequences $\{\Lambda_k\}_{k=1}^\infty, \{\Gamma_k\}_{k=1}^\infty$ of finite sets such that
        $$\frac{105}{64}\delta^k\le\sum_{n\in \Lambda_k}a_n\sign(b_n)\le \frac{153}{64} \delta^k,        \qquad \frac 78\delta^k \le\sum_{n\in \Lambda_k}|b_n|\le \frac 98 \delta^k$$
    and
        $$\frac 78\delta^k \le\sum_{n\in \Gamma_k}|\alpha_n|\le \frac98 \delta^k, \qquad \frac{161}{64}\delta^k\le\sum_{n\in \Gamma_k}\beta_n\sign(\alpha_n)\le \frac{225}{64} \delta^k. $$
 \end{lem}
 \proof
    Using Lemma \ref{lem3.1} for $\{c_n=a_n+ib_n\}_{n=1}^\infty$ and taking $\eta_k=\frac 18\delta^k$, we have
        $$\frac 78\delta^k \le\sum_{n\in \Lambda_k}|b_n|\le \frac 98 \delta^k $$
    and
        $$(2-\frac 18)\frac 78\delta^k\le (2-\frac 18\delta^k)\frac 78\delta^k\le\sum_{n\in \Lambda_k}a_n\sign(b_n)\le (2+\frac18\delta^k)\frac 98 \delta^k\le (2+\frac 18)\frac 98\delta^k.$$
    Similarly, using Lemma \ref{lem3.1} for $\{\beta_n+i\alpha_n\}_{n=1}^\infty$ and taking $\eta_k=\frac 18\delta^k$, we have
        $$\frac 78\delta^k \le\sum_{n\in \Gamma_k}|\alpha_n|\le \frac 98 \delta^k$$
    and
        $$(3-\frac 18)\frac 78\delta^k\le\sum_{n\in \Gamma_k}\beta_n\sign(\alpha_n)\le  (3+\frac 18)\frac 98\delta^k.\eqeproof$$
\bigskip

Now we construct the desired Moran function system from Lemma \ref{lem3.2}. Let $\{\Lambda_k\}$ and $\{\Gamma_k\}$ be given in Lemma \ref{lem3.2}. We define
    $$a^1_k=\sum_{n\in \Lambda_k}a_n\sign(b_n), \qquad b^1_k=\sum_{n\in \Lambda_k}|b_n|$$
and
    $$\alpha^1_k=\sum_{n\in \Gamma_k}|\alpha_n|, \qquad\beta^1_k=\sum_{n\in \Gamma_k}\beta_n\sign(\alpha_n)$$
for $k\ge 1$.

Write
    $$\D_k=\{\pm[(a^1_k+ib_k^1)\pm(\alpha_k^1+i\beta_k^1)]\} =\{d_{k,\,1},d_{k,\,2},d_{k,\,3},d_{k,\,4}\}$$
and then define
    $$f_{k,d}(z)=\delta z+\delta^{1-k}d, \qquad d\in\D_k.$$
We will show that the above function sequence satisfies our demand. Let $Q=[-5, 5]\times [-5, 5]$. It is easy to see that, for all $k\ge 1$,
    \begin{equation}\label{ine1}
        Q\subseteq\bigcup_{d\in\D_k}f_{k, d}(Q).
    \end{equation}

In fact, by Lemma \ref{lem3.2} we have $(a^1_k+\alpha^1_k,
b_k^1+\beta^1_k)\subset [0, 5\delta^k]\times [0, 5\delta^k]$ and
$(a^1_k-\alpha^1_k, b_k^1-\beta^1_k)\subset [0,
5\delta^k]\times[-5\delta^k, 0].$ Then we have $[0, 5]\times [0,
5]\subseteq f_{k,d}(Q)$ if
$d=(a^1_k+ib_k^1)+(\alpha_k^1+i\beta_k^1)$ and $[0, 5]\times [-5,
0]\subseteq f_{k, d}(Q)$ if
$d=(a^1_k+ib_k^1)-(\alpha_k^1+i\beta_k^1)$. This implies
\eqref{ine1} by symmetry property of all $\D_k$.

\begin{lem} \label{lem3.3}
    Let $\{c_n=a_n+ib_n\}_{n=1}^\infty\not\in
    {\ell}^1$ and $\{\gamma_n=\alpha_n+i\beta_n\}_{n=1}^\infty\not\in
   {\ell}^1$ with
    $\{c_n\},\{\gamma_n\}=o(1)$. Suppose $\lim_{n\to\infty}a_n/b_n=2$ and $\lim_{n\to\infty}{\beta_n/\alpha_n}=3$, then there
   exists $c=a+bi\in
   \CC$ such that
        $$c+[-5, 5]\times [-5, 5]\subseteq R(\{c_n\})+R(\{\gamma_n\}).$$
\end{lem}
\begin{proof}
    With the same notations we have a Moran function system $G:=\{f_{k, \, d}(x)=\delta z+\delta^{1-k}d: d\in\D_k, k\ge 1\}$. Note that $f_{k, \, d}(0)=\delta^{1-k}d\in [-5, 5]\times [-5, 5]$. By Proposition \ref{prop3.1} and \eqref{ine1} we have
        $$Q=[-5, 5]\times [-5, 5]\subseteq K_G=\bigg \{f_\sigma(0): \sigma\in\prod_{k=1}^\infty\{1,2,\ldots,n_k\} \bigg\}.$$
    Since
        $$f_{\sigma_1\cdots\sigma_m}(0)=d_{1,\, \sigma_1}+ d_{2,\,\sigma_2}+\cdots+d_{k,\,\sigma_m}, $$
    we have
        $$K_G=\bigg\{\sum_{k=1}^\infty d_{k,\,\sigma_k}: \mbox{all}\,\,\sigma_k\in\{1,2,3,4\}\bigg\} \subseteq R(\{c_n\}_{n\in\cup_{k=1}^\infty\Lambda_k}) +R(\{\gamma_n\}_{n\in\cup_{k=1}^\infty\Gamma_k}).$$
    According to Proposition \ref{coro2.1}, there exist $x_n\in\{-1,1\}$ for $n\in\N \setminus\cup_{k=1}^\infty\Lambda_k$ and $y_n\in\{-1, 1\}$ for $n\in \N\setminus\cup_{k=1}^\infty\Gamma_k$ such that
        $$\sum_{n\in \N\setminus\cup_{k=1}^\infty\Lambda_k}x_nc_n +\sum_{n\in\N\setminus \cup_{k=1}^\infty\Gamma_k}y_n\gamma_n=a+bi.$$
    All the above information together implies the assertion.
\end{proof}

\bigskip

{\noindent {\bf Proof of Theorem \ref{thm2}.}}
    We first note that the sequence $\{c_n\}$ with two ratios must be linearly non-summable.

    By Proposition \ref{prop2.1} we can assume that both 2 and 3 are ratios. Since any  sequence can be decomposed into two sequence with the same ratio,
     we can assume that the sequence$\{c_n\}_{n=1}^\infty$ is decomposed into three sequences: $\{c^{(1)}_n=a^{(1)}_n+ib^{(1)}_n\}_{n=1}^\infty$, $\{c^{(2)}_n=a^{(2)}_n+ib^{(2)}_n\}_{n=1}^\infty$ and $\{c^{(3)}_n\}_{n=1}^\infty$, which satisfy that
        $$\lim_{n\to\infty}a^{(1)}_n/b^{(1)}_n=2, \qquad\lim_{n\to\infty}a^{(2)}_n/b^{(2)}_n=3$$
    and $\{c^{(3)}_n\}_{n=1}^\infty$ is linearly non-summable. The result follows by Lemma \ref{lem3.3} and Theorem \ref{thm1}.
 \eproof

\bigskip

\begin{defi}\label{defi1}
   We call a sequence $\{c_n=a_n+ib_n\}_{n=1}^\infty$ with $a_n/b_n\to t$  {\it changeable } if there
    exists  a partition $\{\Lambda_k\}_{k=1}^\infty$ of $\N$   and a sequence $\{x_n\}_{n=1}^\infty\in\{-1, 1\}^\infty$ such that
      the new sequence $\{\sum_{n\in\Lambda_k}x_nc_n\}_{k=1}^\infty\not\in
      {\ell}^1$ and
       $$\lim_{k\to\infty}\frac{\sum_{n\in\Lambda_k}x_na_n} {\sum_{n\in\Lambda_k}x_nb_n}\ne t.$$
\end{defi}
\begin{theo}\label{thm3.6}
    Let $\{c_n=a_n+ib_n\}_{n=1}^\infty$ be linearly non-summable with $a_n/b_n\to t$. Suppose that the sequence $\{c_n\}$ can be decomposed
     into two non-summable sequences such that one of them is changeable, then $R(\{c_n\})={\CC}$.
\end{theo}
\proof
    By Definition \ref{defi1} we can change the sequence $\{c_n\}$ so that it has at least two distinct radios. Then the result follows by Theorem
    \ref{thm2}.
\eproof

\section{Some examples }
\setcounter{equation}{0}

 We have showed that any complex sequence $\{c_n\}_{n=1}^\infty\not\in
 \ell^1$ has at lease one ratio (Proposition \ref{lem2.1}) and
 $R(\{c_n\})=\CC$ if it has two different ratios (Theorem
 \ref{thm2}). The following example say that there exist complex
 sequences with only one ratio which range is the complex space.
\begin{exam}\label{exam1}
    Let $\{c_n=\frac{(-1)^n}{n\ln (n+1)}+\frac in\}_{n=1}^\infty$ with one ratio $0$. Then $R(\{c_n\})=\CC$.
\end{exam}
\proof
    We decompose the sequence $\{c_n\}$ into two sequences: $\{c_{4k+1}, c_{4k+2}\}_{k=0}^\infty$ and $\{c_{4k+3}, c_{4k+4}\}_{k=0}^\infty$. Note that
        $$-c_{4k+1}+c_{4k+2}=\frac 1{(4k+1)\ln (4k+2)}+\frac 1{(4k+2)\ln(4k+3)}-\frac i{(4k+1)(4k+2)}: =\alpha_k+i\beta_k$$
    and
        $$\lim_{k\to\infty}\frac{\alpha_k}{\beta_k}=\infty.$$
    Then the sequence $\{c_{4k+1}, c_{4k+2}\}_{k=0}^\infty$ is changeable and thus the assertion follows by Theorem \ref{thm3.6}.
\eproof
\begin{rem}
    Let $\{c_n=\frac 1{n\ln (n+1)}+\frac in\}_{n=1}^\infty$ with one ratio $0$. It is easy to show that $R(\{c_n\})$ is dense in $\CC$, but we do not know whether $R(\{c_n\})=\CC$.
\end{rem}

To give an example so that the range of a sequence is dense in $\CC$ but not equal to $\CC$, we begin with the following lemma.

\begin{lem}\label{lem4.1}
    Let $A_n=2^{-n}(\Z+[-1/4, 1/4])$ and $A=\cup_{n=1}^\infty A_{n}$. Then $A\ne \R$.
\end{lem}

\proof
    We show that $1/3\not\in A$. If $1/3\in A$, there exists $n\ge 1$ such that $1/3\in A_n$, i.e., $2^n/3\in \Z+[-1/4, 1/4]$. Note that $2^n/3=k_n+r_n/3$ for some $k_n\in \Z$ and $r_n\in\{1,2\}$, then $\mbox{dist}(2^n/3, \Z)=1/3$. This yields a contradiction to $2^n/3\in \Z+[-1/4, 1/4]$.
\eproof

\begin{exam}\label{exam2}
    Let $\{m_k\}_{k=0}^\infty$ and $\{n_k\}_{k=0}^\infty$ be two increasing integer sequences with $m_0=n_0=0$ and $n_{k+1}\ge n_k+m_k+3$ for $k\ge 0$.
    Define
        $$a_j=2^{-m_k}, \quad b_j=2^{-m_k-n_k} \quad \mbox{and}\, \, c_j=a_j+ib_j,$$
    whenever
        $$\sum_{l=0}^{k-1}2^{m_l+n_l}\le j <\sum_{l=0}^{k}2^{m_l+n_l}\quad \mbox{and}\, \, k\ge 1.$$
    Then $R(\{c_j\})$ is dense in ${\CC}$ but not equal to ${\CC}$.
\end{exam}
\proof
    We first show that the complex sequence $\{c_j\}$ is linearly non-summable.  Note that
        $$\sum_{j=1}^\infty a_j\ge\sum_{j=1}^\infty b_j= \sum_{k=1}^\infty \sum_{\sum_{l=0}^{k-1}2^{m_l+n_l}} ^{\sum_{l=0}^{k}2^{m_l+n_l}-1}
        2^{-m_k-n_k} =\sum_{k=1}^\infty1=\infty.$$
    For any $\alpha, \beta\in\R$ with $\alpha+i\beta\ne0$, clearly $\{\alpha a_j+\beta b_j\}$ does not lie in
    ${\ell}^1$ when $\alpha=0$. When $\alpha\ne 0$, we have $\lim_{j\to\infty}(\alpha a_j+\beta b_j)/a_j=\alpha$. Then $\{\alpha a_j+\beta b_j\}$
    is non-summable and thus $\{c_j\}$ is linearly non-summable. By Theorem \ref{thm1}, $R(\{c_j\})$ is dense in $\CC$.

    Secondly, we show that $R(\{c_j\})$ is not equal to $\CC$. For any $\{x_j\}_{j=1}^\infty\in \{-1, 1\}^\infty$ such that
    $\sum_{j=1}^\infty x_jc_j$ converges, we have
        $$\sum_{j=1}^\infty x_jc_j=\sum_{k=1}^\infty l_k2^{-m_k}+i\sum_{k=1}^\infty l_k2^{-m_k-n_k}, $$
    where all $l_k$ are integers. Hence, there exists $k_0\ge 1$ such that $|l_k2^{-m_k}|\le 1$ for $k>k_0$. Since
        $$\sum_{k=1}^\infty l_k2^{-m_k-n_k}= 2^{-m_{k_0}-n_{k_0}}
        \bigg(
            \sum_{k=1}^{k_0} l_k2^{m_{k_0}-m_k+n_{k_0}-n_k}+ 2^{m_{k_0}}\sum_{k=k_0+1}^\infty l_k2^{-m_k}2^{-(n_k-n_{k_0})}
        \bigg)$$
    and
        $$\bigg|2^{m_{k_0}}\sum_{k=k_0+1}^\infty l_k2^{-m_k}2^{-(n_k-n_{k_0})}\bigg|\le 2^{-(n_{k_0+1}-n_{k_0}-m_{k_0}-1)}\le
       \frac14, $$
    we have $\sum_{k=1}^\infty l_k2^{-m_k-n_k}\in A$, where $A$ is given in Lemma \ref{lem4.1}. Consequently, the imaginary part of $R(\{c_j\})$ is contained in $A$ and thus the assertion follows by Lemma \ref{lem4.1}.\eproof

\section{Hausdorff dimension of the level sets }
\setcounter{equation}{0}

Let $\{x_n\}, \{y_n\}\in\{-1, 1\}^\N$. Define
    $$d(\{x_n\}, \{y_n\})=2^{-k}, $$
where $k$ satisfies that $x_i=y_i$ for $1\le i<k$ and $x_k\ne y_k$.
It is well-known (easy check) that $\{-1, 1\}^\N$ is a complete
metric space with this metric $d(\cdot, \cdot)$. Similarly we define
the Hausdorff dimension on $\{-1, 1\}^\N$ by, for any $B\subseteq
\{-1, 1\}^\N$,
    $$\dim_HB=\sup
    \bigg\{
        s: \lim_{\delta\to 0}\inf
        \Big\{
            \sum_{i\in I}\mbox{diam}(U_i)^s: \{U_i\}_{i\in I} \, \mbox{is a cover of $B$ with $\mbox{diam}(U_i)<\delta$}
        \Big\}=\infty
    \Bigg\}.$$
We begin with a generalization of Theorem 2 in \cite{xi}.
\begin{lem}\label{lem5.1}
    Let $\{c_n\}_{n=1}^\infty$ be a complex number sequence with $\{c_n\}=o(1)$. Then
        $$\dim_H\bigg\{\{x_n\}\in\{-1, 1\}^\N: \sum_{n=1}^\infty x_nc_n\, \mbox{converges}\bigg\}=1.$$
\end{lem}
\proof
    It is easy to see that the proof of Theorem 2 in \cite{xi} will be used here by appropriate modifications.
\eproof

\begin{defi}
    Let $\Lambda\subseteq \N$. Define the {\it super and lower density} of $\Lambda$ by
        $$\overline{D}(\Lambda) =\limsup_{k\to\infty}\frac{\#(\Lambda\cap[0, k])}{k}, \qquad \underline{D}(\Lambda)= \liminf_{k\to\infty}\frac{\#(\Lambda\cap[0, k])}{k}, $$
    respectively, where $\#E$ is the cardinalities of the set $E$. If $\overline{D}(\Lambda)=\underline{D}(\Lambda)$, we say the common value the {\it density} of $\Lambda$ and denote it by $D(\Lambda)$.
\end{defi}

Let $\Lambda=\{n_1, n_2, \cdots\}\subseteq\N$. We define a map $h_\Lambda$ from $\{-1, 1\}^\N$ to itself by
    $$h_\Lambda(\{x_n\}_{n\in\N})=\{x_n\}_{n\in \N\setminus \Lambda}.$$

\begin{lem} \label{lem5.2}
    Let $\Lambda\subseteq\N$ and $0<\epsilon<1$. If $\overline{D}(\Lambda)<\epsilon$, then
        $$(1-\epsilon)\dim_Hh_\Lambda(B)\le\dim_H B$$
    for any $B\subseteq \{-1, 1\}^\N$.
\end{lem}
\proof
    Denote $m_k=\#(\Lambda\cap [0, k-1])$. Then, by $\overline{D}(\Lambda)<\epsilon$, $m_k/k<\epsilon$ for $k>k_0\ge 1$. For any $\{x_n\}, \{y_n\}\in \{-1, 1\}^\N$  with $d(\{x_n\},\{y_n\})=2^{-k}<2^{k_0}$, we have
        $$d\big(h_\Lambda(\{x_n\}), h_\Lambda(\{y_n\})\big)\le 2^{-k+m_k}=d(\{x_n\}, \{y_n\})^{1-\frac{m_k}k}<d(\{x_n\}, \{y_n\})^{1-\epsilon}.$$
    Hence, by the definition of Hausdorff dimension, it is easy to check that $(1-\epsilon)\dim_Hh_\Lambda(B)\le\dim_H B $.
\eproof

\begin{lem}
    Let $\{c_n=a_n+ib_n\}_{n=1}^\infty$ be a linearly non-summable sequence and let $\epsilon$
    so that $0<\epsilon<1$. Suppose that $\{c_n\}$ has a unique ratio (at least two distinct ratios), then there exists a linearly non-summable subsequence $\{c_n\}_{n\in\Lambda}$ with one ratio (two distinct ratios, resp.) so that $\overline{D}(\Lambda)<\epsilon$.
\end{lem}
\proof
    First we show the case of one ratio. By Proposition \ref{prop2.1} we can assume that the unique ratio is $1$, that is, $\lim_{n\to\infty}a_n/b_n=1$. Since $\sum_{n=1}^\infty |a_n-b_n|=\sum_{j=1}^q\sum_{k=0}^\infty|a_{kq+j}-b_{kq+j}| =\infty$,
    there exists $j$ so that $\sum_{k=0}^\infty|a_{kq+j}-b_{kq+j}|=\infty$. This implies that the subsequence $\{c_{kq+j}\}_{k=0}^\infty$ is linearly non-summable. Denote $\Lambda_q=\{kq+j: k=0, 1, 2, \ldots\}$, by a simple calculation we have $\overline{D}(\Lambda_q)=1/q$. Hence the assertion follows by choosing $q$ so that $1/q<\epsilon$.
    Secondly, for the case of at least two ratios, the assertion follows by the same idea used for two subsequences with distinct ratios.
\eproof

\begin{theo}
    Let $\{c_n=a_n+ib_n\}_{n=1}^\infty$ be a linearly non-summable sequence and let $c\in\CC$.
    \begin{enumerate}
    \item If $\{c_n\}$ has one ratio, then
            $$\dim_H\bigg\{\{x_n\}\in\{-1, 1\}^\N: \sum_{n=1}^\infty x_nc_n\in B(c, \delta)\bigg\} =1$$
        for any $\delta>0$;
    \item If $\{c_n\}$ has at least two distinct  ratios, then
            $$\dim_H\bigg\{\{x_n\}\in\{-1, 1\}^\N: \sum_{n=1}^\infty x_nc_n=c\bigg\}=1.$$
    \end{enumerate}
\end{theo}

\proof
    (1) For any $\epsilon>0$, by Lemma \ref{lem5.2} there {exists} a linearly non-summable subsequence $\{c_n\}_{n\in\Lambda}$ with  $\overline{D}(\Lambda)<\epsilon$. According to Lemma \ref{lem5.1}, we have
        $$\dim_H\bigg\{\{x_n\}\in\{-1, 1\}^\N: \sum_{n\in\N}x_nc_n\, \mbox{converges}\bigg\}=1.$$
    To show the assertion (1), it is sufficient, by Lemma \ref{lem5.2}, to show that
    \begin{equation}\label{inc1}
        h_\Lambda\bigg(\bigg\{\{x_n\}\in\{-1, 1\}^\N:\sum_{n\in\N}x_nc_n\in B(c, \delta)\bigg\}\bigg)\supseteq \bigg\{\{x_n\}_{n\in\N\setminus\Lambda}: \sum_{n\in\N\setminus\Lambda}x_nc_n\, \mbox{converges}\bigg\}.
    \end{equation}

    For any $\{x_n\}_{n\in\N\setminus\Lambda}$ so that $\sum_{n\in\N\setminus\Lambda}x_nc_n$ converges (to  $d$), by Theorem \ref{thm1} there exists $\{x_n\}_{n\in\Lambda}$ such that $\sum_{n\in\Lambda}x_nc_n\in B(c-d, \delta)$. Then we have
        $$ \sum_{n\in\N}x_nc_n= \sum_{n\in\N\setminus\Lambda}x_nc_n+ \sum_{n\in\Lambda}x_nc_n\in B(c, \delta).$$
    This implies \eqref{inc1} by the definition of $h_\Lambda$.

    The proof of (2) is easy to be given by the similar idea of (1).
 \eproof

\end{document}